\documentclass[]{amsart}
\usepackage{amsmath,amsfonts,amssymb}

\def\C{\mathbb C}
\def\Q{\mathbb Q}
\def\P{\mathbb P}

\def\diff{{\rm Diff}}

\def\Cdos{(\C^2,0)}
\def\ord{\mathop{\rm ord}}

\theoremstyle{plain}
\newtheorem{theorem}{Theorem}[section]
\newtheorem{lemma}[theorem]{Lemma}

\newtheorem{proposition}[theorem]{Proposition}
\def\proof{{\it Proof: }}
\def\qed{\hfill\hbox{$\square$}}

\theoremstyle{definition}

\author[F.E. Brochero Mart\'{\i}nez]{F. E. Brochero Mart\'{\i}nez}
\address{
Departamento de Matem\'atica\\
UFMG\\
Belo Horizonte, MG\\
 30123-970\\
 Brazil\\
 }
 \email{fbrocher@mat.ufmg.br }
\thanks{The first author was  supported by CAPES, Brazil, Process:
BEX3083/05-5}

\author[F. Cano]{F. Cano}

\author[L. L\'opez-Hernanz]{L. L\'opez-Hernanz}
\address{Departamento de \'Algebra, Geometr\'{\i}a y Topolog\'{\i}a\\
Universidad de Valladolid, Spain}
\email{fcano@agt.uva.es}\email{llopez@agt.uva.es}
\thanks{The third author was  supported by FPU, Spain, Process:AP2005-3784}

\date{\today
}

\subjclass[2000]{32H02, 32H50, 37F99}

\title{Parabolic curves for  diffeomorphisms in $(\C^2,0)$}

\begin{document}

\begin{abstract}
We give  a simple proof of the existence of parabolic curves for
tangent to the identity  diffeomorphisms in $\Cdos$ with isolated
fixed point.
\end{abstract}

\maketitle

\section{Introduction}
Let $F$ be a tangent to the identity diffeomorphism of $(\C^n,0)$. A
{\em parabolic curve} for $F$ is an injective holomorphic map
$\varphi:\Omega\to \C^n$, where $\Omega$ is a simply connected
domain in $\C$ with $0\in\partial \Omega$ such that
\begin{enumerate}
\item $\varphi$ is continuous at the origin, and $\varphi(0)=0$.
\item $F(\varphi(\Omega))\subset \varphi(\Omega) $ and  $F^{\circ k}(p)$ converges to $0$  when $k\to
+\infty$, for $p\in \varphi(\Omega)$.
\end{enumerate}
We say that  $\varphi$ is {\em tangent} to $[v]\in \P^{n-1}$  if
$[\varphi(\zeta)]\to [v]$ when $\zeta\to 0$. Let us write
$F(z)=z+P_k(z)+P_{k+1}(z)+\cdots$,   where  $P_j$ is a
$n$-dimensional vector of homogeneous polynomials of degree $j$, and
$P_k(z)\not\equiv 0$.  A {\em characteristic direction} for $F$ is a
point $[v]\in \P^{n-1}$ such that $P_k(v)=\lambda v$, for some
$\lambda\in \C$; it is {\em nondegenerate} if $\lambda\ne 0$. The
integer $\ord(F):=k\ge 2$ is the {\em tangency order} of $F$ at $0$.

The following theorem is  analogous to Briot and Bouquet's theorem
\cite{BrBo} for diffeomorphisms of $(\C^n,0)$.

\begin{theorem}[Hakim \cite{Hak}]
Let $F$ be a  tangent to the identity germ of  diffeomorphism of
$(\C^n,0)$. For any nondegenerate characteristic direction $[v]$
there exist $\ord(F)-1$ disjoint parabolic curves tangent to $[v]$
at the origin.
\end{theorem}
When $n=2$,
Abate proved  that the nondegeneracy condition can be
dismissed.

\begin{theorem}[Abate \cite{Aba,ABT}]\label{abate} Let  $F$ be
a tangent to the identity germ of diffeomorphism of $\Cdos$ such
that $0$ is an isolated fixed point. Then there exist $\ord(F)-1$
disjoint parabolic curves for $F$ at the origin.
\end{theorem}
This theorem is analogous to Camacho-Sad's theorem \cite{CaSa} of
existence of invariant curves for holomorphic vector fields. We show
in this note that the analogy is deeper enough to prove theorem
\ref{abate} in a simple way  starting with Hakim's theorem.

\section{Exponential Operator and Blow-up transformation}

 Let ${\hat{\mathfrak X}}_2\Cdos$  be the module of formal vector
 fields $X=a(x,y)\frac {\partial}{\partial x}+
b(x,y)\frac {\partial}{\partial y}$ of order $\ge 2$, i.e.,
$\min\{\nu(a),\nu(b)\}\ge 2$. We denote by $\widehat\diff_1\Cdos$
the group  of tangent  to the identity formal diffeomorphisms
$F(x,y)=(x+p(x,y),y+q(x,y))$  where
$\min\{\nu(p(x,y)),\nu(q(x,y)\}\ge 2$. Let us denote by ${\mathfrak
X}_2\Cdos$ and by $\diff_1\Cdos$ the convergent elements of
${\hat{\mathfrak X}}_2\Cdos$ and $\widehat\diff_1\Cdos$
respectively.

 Let $X\in {\widehat{\mathfrak X}}_2\Cdos$. The exponential
operator of $X$ is  the application $\exp{tX}:\C[[x,y]]\to
\C[[x,y,t]]$ defined by the formula
$$\exp{tX}(g)=\sum_{j=0}^\infty \frac
{t^j}{j!}X^j(g)$$ where $X^0(g)=g$ and $X^{j+1}(g)=X(X^j(g))$. Note
that, since  $\nu(X^j(g))\ge j+\nu(g)$, we can substitute $t=1$ to
get the element  $\exp{X}(g)\in \C[[x,y]]$. Moreover, $\exp{tX}$
gives a homomorphism of $\C$-algebras, in particular, we have
$$\exp{tX}(fg)=\exp{tX}(f)\exp{tX}(g).$$
We get also
\begin{proposition}\label{expo}The application
$$\begin{array}{rcl}Exp:{\widehat{\mathfrak X}}_2\Cdos&\to&
\widehat\diff_1\Cdos\\
X&\mapsto&(\exp X(x),\exp X(y))
\end{array}$$
is a bijection.
\end{proposition}
\proof Let $G(x,y)=\Bigl(x+\sum\limits_{n=2}^\infty
p_n(x,y),y+\sum\limits_{n=2}^\infty q_n(x,y)\Bigr)$ and $
X=\sum\limits_{n=2}^\infty\Bigl( a_n(x,y)\frac{\partial}{\partial
x}+ b_n(x,y)\frac{\partial}{\partial y}\Bigr)$. The identity
$Exp(X)=G$ is equivalent to
$$\begin{array}{rcl}
p_{m+1}&=&a_{m+1}+HT_{m+1}\Bigl(\sum\limits_{j=2}^{m} \frac
1{j!}X_m^j(x)\Bigr)\\
q_{m+1}&=&b_{m+1}+HT_{m+1}\Bigl(\sum\limits_{j=2}^{m} \frac
1{j!}X_m^j(y)\Bigr),
\end{array}
$$
where $X_m=\sum\limits_{n=2}^m\Bigl(
a_n(x,y)\frac{\partial}{\partial x}+
b_n(x,y)\frac{\partial}{\partial y}\Bigr),$ and $HT_{m+1}(h)$ is the
homogeneous term of $h$ of order $m+1$. These equations determine
univocally $X$ if $G$ is given.\qed

In general, $X$ may be  not convergent for certain convergent $G$. The formal vector field $X$
such that $G=Exp(X)$ is called the {\em infinitesimal generator} of $G$. Note  that $\ord
(G)=\nu(X)$. If $k=\nu(X)$, then  $a_k=p_k$ and $b_k=q_k$, thus  the characteristic directions
of $F$ correspond to the points of  the tangent cone of $X$. Moreover,  if
  $X=fX' $ with $X'\in {\widehat{\mathfrak X}}\Cdos$ and $f\in \C[[x,y]]$
then  $Exp(X)(x,y)=(x+f(x,y)p(x,y),y+f(x,y)q(x,y))$. The converse
statement follows by a  process similar to the proof of proposition
\ref{expo}. In particular, $0$ is an isolated singular point of $X$
if and only if $0$ is an isolated fixed point of $F$. In the case
$f(x,y)=x^k$, and $S=(x=0)$  invariant by $X'$, Camacho-Sad's index
of $X$ at $0$ along $S$ is exactly Abate's residual  index of $F$ at
$0$ along $S$.

Now, let $\pi:(M,D)\to \Cdos$ be the blow up of $\C^2$ at the
origin, where $D=\pi^{-1}(0)=\P^1$, thus  each characteristic
direction determines a point of $D$.

\begin{proposition}
Let $F\in \diff_1\Cdos$. There exists a unique
 germ of diffeomorphism $\tilde F$ in $(M,D)$ such that  $\pi\circ \tilde F=F\circ
\pi$ and  $\tilde F|_D=id|_D$. Moreover, the germ $\tilde F_p$ has
order $\ge \ord(F)$ for any characteristic direction $p\in D$ and
hence $\tilde F_p\in \diff_1(M,p)$.
\end{proposition}
\proof Let $F(x,y)=(x+p_k(x,y)+\cdots,y+q_k(x,y)+\cdots)$ where
$k=\ord(F)\ge 2$. We have   two charts of $M=U_1\cup U_2$ such that
$\pi|_{U_1}:U_1\to \C^2$, is defined by $\pi(x,v)=(x,xv)$ and
$\pi|_{U_2}:U_2\to \C^2$, is defined by $\pi(u,y)=(uy,y)$. We define
$\tilde F$ in the first chart as
$$\begin{array}{rcl}
\tilde F(x,v)&=&\pi^{-1}\circ F\circ
\pi(x,v)=\Bigl(x+p_k(x,xv)+\cdots,\dfrac{vx+q_k(x,xv)+\cdots}{x+p_k(x,xv)+\cdots}\Bigr)\\
&=&(x+x^{k}(p_k(1,v)+x(\cdots)),v+x^{k-1}(q_k(1,v)-vp_k(1,v)+x(\cdots)))
\end{array}
$$
Observe that  $\tilde F(0,v)=(0,v)$, thus any point of the divisor
is fixed. Moreover, if  $q_k(1,v_0)-v_0p_k(1,v_0)=0$ we have
$dF(0,v_0)=I$, and thus for any characteristic direction
$p=(0,v_0)\in D$, $\ord(\tilde F_p)\ge \ord(F)$.\qed

\begin{proposition}
Let $X\in{\widehat{\mathfrak X}}_2\Cdos$. Let $\tilde X$ be the
formal vector field in $(M,D)$ such that $D\pi\cdot\tilde X=X\circ
\pi$. If $p$ is a point of the tangent cone of $X$  then $\tilde
X_p\in{\widehat{\mathfrak X}}_2(M,p)$.
\end{proposition}

\proof  Let $X=a(x,y)\frac{\partial}{\partial
x}+b(x,y)\frac{\partial}{\partial y}$ with $a(x,y)=a_k(x,y)+\cdots$,
$b(x,y)=b_k(x,y)+\cdots$ and $k\ge 2$. Let $U_1$ and $U_2$ be two
charts of $M=U_1\cup U_2$ as in the proposition above. Then $\tilde
X$ is given in the chart $U_1$ by
$$\begin{array}{rcl}\tilde X(x,v)
&=&a(x,xv)\frac{\partial}{\partial x}+\dfrac {b(x,xv)-va(x,xv)}x\frac{\partial}{\partial v}\\
&=&x^k(a_k(1,v)+x(\cdots))\frac{\partial}{\partial
x}+x^{k-1}((b_k(1,v)-va_k(1,v))+x(\cdots))\frac{\partial}{\partial
y}
\end{array}
$$
Now, if $p=(0,v_0)\in D$ is such that $b_k(1,v_0)-v_0a_k(1,v_0)=0$,
then $\nu_p(a(x,xv))\ge k$ and $\nu_p\bigl(\frac
{b(x,xv)-va(x,xv)}x\bigr)\ge k$  so $\tilde
X_p\in{\widehat{\mathfrak X}}_2(M,p)$. \qed

We say that $X$ is {\em strictly non singular} if $X=fX'$, where
$X'$ is a non singular formal vector field. Otherwise, we say that
$X$ is {\em strictly singular}. Note that in the above statement any
strictly  singular point of $\tilde X$ is in the tangent cone of
$X$. Let us also recall that Seidenberg's reduction of singularities
\cite{Sei} is done by blowing-up at strictly singular points.

\begin{lemma}\label{comuta} Let $F\in \diff_1\Cdos$  and $X\in {\widehat{\mathfrak
X}}_2\Cdos$ such that $F=Exp(X)$. Let  $\tilde X$ be as in the
proposition above. Then for any characteristic direction $p\in D$
$$\tilde F_p=Exp (\tilde X_p).$$
\end{lemma}

\proof Let $U\simeq\C^2$ be a chart of $M$ such that $\pi|_U:U\to
\C^2$ is defined by $\pi(x,v)=(x,xv)$ and  $p\in U\cap D=\{(0,v)\in
U\}$ be a point on the divisor.   Without lost of generality,
applying a linear change of coordinates, we can suppose that
$p=(0,0)\in U$. Since
$$F(x,y)=Exp(X)=(\exp X(x),\exp X(y)),$$
using the definition of $\tilde F$, we have
$$\begin{array}{rcl}\tilde F(x,v)&=&\Bigl(\exp X(x), \dfrac {\exp X(xv)}{\exp X(x)}\Bigr)=
\Bigl(\exp X(x), \dfrac {\exp X(x)\exp X(v)}{\exp X(x)}\Bigr)\\
&=&(\exp X(x),\exp X(v))=(\exp {\tilde X}(x),\exp {\tilde X}(v))
=Exp(\tilde X_p)(x,v).
\end{array}
$$
This ends the proof.\qed

\section{Existence of parabolic curves}

We need the following formal version of Camacho-Sad's theorem
\cite{CaSa} whose proof goes exactly as the original one (see also
\cite{CanJ}).

\begin{theorem}[Camacho and Sad] Take  $X\in {\widehat{\mathfrak
X}}_1\Cdos$  with an isolated singularity at the origin. There is a
desingularization morphism $\sigma:(\tilde M, \tilde D)\to \Cdos$
composition of a finite sequence of blow-ups with centers at
strictly singular points and a point $ p\in \tilde D$ satisfying the
following property: There are local coordinates $(u,v)$ at $p$ such
that $\tilde D_p=(u=0)$ and the transform $X^*$  of $X$ at $p$ is of
the form:
$$X^*(u,v)=u^k\Bigl((\lambda u+u^2(\cdots))\frac{\partial }{\partial
u}+(\mu v+ u(\cdots))\frac{\partial }{\partial v}\Bigr)$$ where
$\lambda\ne 0$ and $\frac {\mu}{\lambda}\notin \Q_{>0}$.
\end{theorem}

Let us prove theorem \ref{abate}. Take $X$  the infinitesimal
generator of $F$, and consider  $X^*$ and $p$ as in Camacho-Sad's
Theorem. By lemma \ref{comuta} we have
$$F^*_p(u,v):=Exp(X^*_p)=(u+\lambda u^{k+1}+O(u^{k+2}), v+\mu
u^kv+O(u^{k+1}))$$ so $F^*_p$ is a diffeomorphism tangent to the
identity, with $(1,0)$ as a nondegenerate characteristic direction.
By  Hakim's Theorem, there exist $\ord(F)-1$ disjoint parabolic
curves $\varphi_j:\Omega_j\to \tilde M$ for $F^*_p$ tangent to the
direction $(1,0)$ at $p$. Since this direction is transversal to the
divisor, it follows that $\overline{\varphi_j(\Omega_j)}\cap \tilde
D=\{p\}$ and thereby $\pi\circ\varphi_j$ is also a parabolic curve
for $F$.

Furthermore, according to J. Cano's proof \cite{CanJ} of
Camacho-Sad's theorem, to find  the points $p\in \tilde D$ that
satisfy Camacho-Sad's theorem, it is enough to follow after the
first blow up, the singularities with Camacho-Sad's index not in
$\Q_{\ge 0}$. Thus, there exist parabolic curves for any
characteristic direction of $F$ that gives at the divisor Abate's
residual index not in $\Q_{\ge0}$ (see corollary 3.1. in
\cite{Aba}).


\begin{thebibliography}{99}
\bibitem{Aba} Abate, M., {\it The residual index and the dynamics of
holomorphic maps  tangent to the identity}, Duke Math. J. {\bf 107}
(2001) 173-207

\bibitem{ABT}
 Abate M., Bracci F., Tovena F., {\it Index theorems for holomorphic
self-maps}. Ann. Math. {\bf 159} (2004), 819-864


\bibitem{BrBo} Briot, Bouquet, {\it Recherches sur les propri\'et\'es les fonctions
d\'efinies par des \'equations diff\'erentielles}, J. Ecole
Polytechnique, XXI, (1856) 133-198

\bibitem{CaSa} Camacho, C., Sad P., {\it Invariant varieties through singularities of holomorphic vector fields} Ann. of Math. {\bf 115}
(1982) 579-595

\bibitem{CanJ} Cano, J., {\it Construction of invariant curves for
singular holomorphic vector fields}, Proc. Amer. Math. Soc., {\bf
125} (1997) 2649-2650

\bibitem{Hak} Hakim, M., {\it Analytic transformations of $(\C^p,0)$
tangent to the identity}, Duke Math. J. {\bf 92} (1998) 403-428.

\bibitem{Sei} Seidenberg, A., {\it Reduction of singularities of the
differential equation  $AdY=BdX$}, Amer. J. Math. (1968) 248-269





\end{thebibliography}
\end{document}